\documentclass[review]{elsarticle}
\usepackage{lineno,hyperref}
\modulolinenumbers[5]
\journal{}

\usepackage{multirow, booktabs}
\usepackage{amsmath,amsthm}
\usepackage{amssymb,latexsym}
\usepackage{graphicx}
\usepackage[labelsep=endash]{caption}
\usepackage{mathptmx}

\usepackage{grffile}
\usepackage{relsize}
\usepackage{multirow}
\usepackage{subfig}
\newtheorem{thm}{Theorem}[section]
\newtheorem{cor}{Corollary}[section]
\newtheorem{lem}{Lemma}[section]
\newtheorem{prop}{Proposition}[section]
\theoremstyle{definition}

\newtheorem{rem}{Remark}[section]

\begin{document}
\begin{frontmatter}
\title{New wavelet method based on Shifted Lucas polynomials: A tau approach}

\author[label1]{Rakesh Kumar*}
\ead[a]{rakesh.lect@gmail.com}
\author[label1]{Reena Koundal}
\author[label1]{K. Srivastava}
\ead[b]{Supported by DST, Ministry of Science and Technology, India through WOS-A vide their File No.SR/WOS-A/PM-20/2018}
\address[label1]{School of Mathematics, Computer \& Information Sciences, Central University of Himachal Pradesh, Dharamshala, India}
\begin{abstract} 
In current work, non-familiar shifted Lucas polynomials are introduced. We have constructed a computational wavelet technique for solution of initial/boundary value second order differential equations. For this numerical scheme, we have developed weight function and Rodrigues' formula for Lucas polynomials. Further, Lucas polynomials and their properties are used to propose shifted Lucas polynomials and then utilization of shifted Lucas polynomials provides us shifted Lucas wavelet. We furnished the operational matrix of differentiation and the product operational matrix of the shifted Lucas wavelets. Moreover, convergence and error analysis ensure accuracy of the proposed method. Illustrative examples show that the present method is numerically fruitful, effective and convenient for solving differential equations
\end{abstract}
\end{frontmatter}
\section{Introduction}
 The various properties like orthogonality, compact support, arbitrary regularity and high order vanishing moments make the theory of wavelets more powerful. With all of these properties, wavelets have attracted the consideration of every researcher due to their applications in wave propagation, pattern recognition, computer graphics and medical image technology \cite{Constantmldes1987}. The quality smoothness and better interpolation are manufacturing the wavelets technique based on orthogonal polynomials, more friendly to get the solutions of different kinds of differential equations. We can also say that other than numerical methods, wavelets theory provide us a new direction for solutions of the differential equations.\\
In the present work, we are interested in the solutions of two different kinds of second order differential equations. 
Our main motivation behind the present work is to construct a new technique based on orthogonal Lucas polynomials and wavelets transformation for solutions of the Lane-Emden and Pantograph differential equations. To the best of our knowledge, no one has introduced Lucas wavelets yet. \\
The remaining structure is as follows: section 2, devoted to an overview on Lucas polynomials. In section 3, we discussed about shifted Lucas polynomials and their orthogonality condition. In section 4, wavelets and their properties are presented. In section 5, procedure for Lane-Emden and Pantograph differential equation is discussed. Convergence and Error analysis are derived in section 6. The validity of our method is checked by numerical experiment in section 7. Conclusion is done in the last section.
\section{Overview of Lucas polynomials}
\noindent A. F. Horadam \cite{AFH} presented the polynomial sequence $\{W_{s}(\Theta)\}_{s\in\mathbb{N}\cup\{0\}}$ defined by the 
relation 
$$W_s(\Theta)=\vartheta(\Theta)W_{s-1}(\Theta)+\upsilon(\Theta)W_{s-2}(\Theta),\ s\in\{2, 3, 4,\ldots\},$$ with initial conditions: $W_0(\Theta)=\chi_0,W_1(\Theta)=\chi_1\Theta^K$, where $ \vartheta(\Theta)=\chi_2\Theta^K,\ \upsilon(\Theta)= \chi_3\Theta^K$ in which $\chi_0, \chi_1, \chi_2, \chi_3$ are constants and $K=0$ or $1$. If $W_0(\Theta)=2$ and $W_1(\Theta)=\vartheta(\Theta)$ then $W_s(\Theta)$ is expressed as $W_s(\Theta)=\chi^s(\Theta)+\wp^s(\Theta)$, where $\chi(\Theta)+\wp(\Theta)=\vartheta(\Theta)$ and $\chi(\Theta)\wp(\Theta)=-\upsilon(\Theta)$. Here, $\{W_s\}$ is recognized as Lucas polynomial sequences. The Lucas polynomials are defined by adopting $\vartheta(\Theta)=\Theta$ and $\upsilon(\Theta)=1$ in $\{W_s\}$. It is given explicitly as follows: 
\begin{align}\label{lucas}
L_s^{\ast}(\Theta)=2^{-s}\left[\left( \Theta-\sqrt{\Theta^2+4}\right)^{s} +\left( \Theta+\sqrt{\Theta^2+4}\right)^s\right].
\end{align}
Alternatively, the Lucas polynomials $ L_{s}^{*}(\Theta) $ of degree $ s $ can also be defined by the following relation \cite{Lachal2013}: 

\begin{align}\label{1.1}
L_{s}^{*}(2\sinh\theta)=
\begin{cases}
2\sinh(s\theta),\hspace{0.2cm} \text{if $s$ is odd, }\\
2\cosh(s\theta),\hspace{0.2cm} \text{if $s$ is even.}                
\end{cases}
\end{align}
Some of the Lucas polynomials are as follows: $L_{0}^{*}(\Theta)=2$, $L_{1}^{*}(\Theta)=\Theta$, $L_{2}^{*}(\Theta)=\Theta^{2}+2$,  $L_{3}^{*}(\Theta)=\Theta^{3}+3\Theta$, $L_{4}^{*}(\Theta)=\Theta^{4}+4\Theta^{2}+2$ and so on.

\begin{prop} The differential equation 
\begin{align}\label{1.2}
(\Theta^{2}+4)\dfrac{d^{2}y}{d\Theta^{2}}+\Theta\dfrac{dy}{d\Theta}-s^{2}y=0
\end{align}
is satisfied by the Lucas polynomial.
\end{prop}

\begin{proof}
Proof is omitted.
\end{proof}

\noindent The differential equation given by Eq. \eqref{1.2} is called the {\it Lucas differential equation}. From the above discussion, it is not hard to see the following proposition:

\begin{prop}
If $L_{0}^{*}(\Theta)=2$ and $L_{1}^{*}(\Theta)=\Theta$, then three-term recurrence relation of Lucas polynomial  $ L_{s}^{*}$ is defined as
\begin{align}\label{rec}
L_{s}^{*}(\Theta) = \Theta L_{s-1}^{*}(\Theta)+L_{s-2}^{*}(\Theta),\ s\in\{2,3,4,\ldots\},
\end{align}
\end{prop}

\noindent A sequence of polynomial $\{f_s(\Theta)\}$ is said to be orthogonal on a rectifiable Jordan curve $\Gamma$ if the following relation holds:
$$\dfrac{1}{\ell}\int_{\Gamma}f_s(\Theta)\overline{f_h(\Theta)}\omega(\Theta)\vert d\Theta\vert=\delta_{sh},$$
where $w(\Theta)$ is a weight function defined on $\Gamma$ and $\ell$ is the perimeter of the boundary of $\Gamma$ in the complex plane \cite{GS}.

\begin{lem}\label{lemma}
The Lucas polynomials $ L_{s}^{*}(\Theta)$ are orthogonal with respect to weight function $ w(\Theta) $.
\end{lem}
\begin{proof}
\begin{align*}
\langle L_{s}^{*}, \overline{L_{h}^{*}}\rangle = \begin{cases}
0,\hspace{0.6cm} \text{if $s\neq h$, }\\
\pi/2,\hspace{0.2cm} \text{if $s=h>0$, }\\
\pi,\hspace{0.53cm} \text{if $s=h=0$.}                
\end{cases} \end{align*}
\end{proof}
\begin{rem}
The Lucas polynomials $ L_{s}^{*}(\Theta) $ have precisely $s$ zeroes in the form \citep{Koshy2017} 
\begin{align*}
\Theta_{j} = \left(2\iota\cos\dfrac{(2j+1)\pi}{2s} \right); j=0,1,2,\ldots, s-1.
\end{align*} 
\end{rem}
\begin{thm}
\textbf{(Rodrigues' Formula)}. The Lucas polynomials $L_{s}^{*}(\Theta)$ can be acquired through the Rodrigues' formula given by:\\
\begin{align}\label{Rodrigues'}
L_{s}^{*}(\Theta)=2\dfrac{s!}{(2s)!}\left({\Theta^{2}+4}\right)^{1/2}\dfrac{d^{s}}{d\Theta^{s}}\left\{\left({\Theta^{2}+4}\right)^{s-1/2}\right\}. 
\end{align}
\end{thm}
\begin{proof}
Proof is omitted.
\end{proof}

\begin{prop}
The generating function for Lucas polynomial is defined as \cite{Boussayoud2017}
\begin{align*}
\sum_{s=0}^{\infty}L_{s}^{*}(\Theta)t^{h}= \dfrac{2-\Theta t}{1-\Theta t-t^{2}}.
\end{align*}
\end{prop}
\section{Shifted Lucas polynomials}
The shifted Lucas polynomials are equally important as Lucas polynomials. Few shifted Lucas polynomials are given below:
\begin{align*}
\begin{aligned}
&Q_{0}^{*}(t)=2,&\\&
Q_{1}^{*}(t)=2t-2\iota,&\\
&Q_{2}^{*}(t)=4t^{2}-8\iota t-2,&\\
&Q_{3}^{*}(t)=8t^{3}-24\iota t^{2} t-18t+2\iota ,&\\
&\vdots
\end{aligned}
\end{align*}
Also it can be easily seen that the orthogonality condition for the shifted Lucas polynomials is as follows:
\begin{align*}
\langle Q_{s}^{*}, \overline{Q_{h}^{*}}\rangle
=\begin{cases}
0,\hspace{0.6cm} \text{if $s\neq h$, }\\
\frac{\pi \alpha_{s}}{2},\hspace{0.3cm} \text{if $s=h$, }               
\end{cases}
\end{align*}
where
\begin{align}\label{alpha}
\alpha_{s}=\begin{cases}
2,\hspace{0.6cm} \text{if $s=0$, }\\
1,\hspace{0.6cm} \text{if $s\geq 1$.}               
\end{cases}
\end{align}

\section{Lucas wavelets and shifted Lucas wavelet}
Analogous to Legendre wavelets, Chebyshev wavelets , Lucas wavelets also have  $\varphi_{h,s}(\Theta)= \varphi(k,h,s,\Theta)$ four arguments. 
On the interval $[0,2)$, the Lucas wavelets family of functions can be defined as
 \begin{align}\label{3}
 \varphi_{h,s}(\Theta)=
\begin{cases}
 2^{\frac{k+1}{2}} \widetilde{{L}_{s}^{*}}(2^{k+1}\iota \Theta-\hat{h}\iota), \hspace{0.6cm} \dfrac{\hat{h}-2}{2^{k+1}}\leq \Theta<\dfrac{\hat{h}+2}{2^{k+1}},\\
0, \hspace{3.4cm} o.w,\\                      
\end{cases}
\end{align}  
where
\begin{align}\label{31}
\widetilde{{L}_{s}^{*}}(\Theta)=
\begin{cases}
\dfrac{1}{\sqrt{\pi}},\hspace{1.2cm}s=0,\\
\sqrt{\dfrac{2}{\pi}}L_{s}^{*}(\Theta), \hspace{0.35cm} s>0,\\                      
\end{cases}
\end{align}
$s=0,1,2,\ldots,S-1;\ S$ is the maximum order of the Lucas polynomials; $\hat{h}=2(2h+1)$ where $ h=0,1,2,\ldots,2^{k}-1 $; $k$ is any nonnegative integer and $L_{s}^{*}(\Theta) $ are Lucas polynomials of order $s$. Also, On the interval $ [0,2] $ shifted Lucas wavelet is defined as
\begin{align*}
\varphi_{h,s}(\Theta)=
\begin{cases}
 2^{\frac{k+1}{2}} \sqrt{\dfrac{2}{\pi\alpha_{s}}} Q_{s}^{*}(2^{k}\iota \Theta-2h\iota), \hspace{0.6cm} \dfrac{h}{2^{k-1}}\leq\Theta<\dfrac{h+1}{2^{k-1}},\\
0, \hspace{3.4cm} o.w,\\                      
\end{cases}
\end{align*}
where $ h=0,1,2,..., 2^{k}-1 $ and $ \alpha_{s} $ is defined in previous definitions.

\section{Shifted Lucas tau wavelet technique for Lane-Emden and Pantograph differential equations}\label{g}
A function $ \rho(\Theta)\in L_{w_{h}}^{2}[0,2] $ can be defined as a Lucas wavelet series in the following form
\begin{align}\label{4}
\rho(\Theta)= \sum_{h=0}^{\infty}\sum_{s=0}^{\infty}\mathcal{E}_{h,s}\varphi_{h,s}(\Theta),
\end{align}
where
\begin{align}\label{5}
\mathcal{E}_{h,s}=\left\langle \rho(\Theta), \varphi_{h,s}(\Theta)\right\rangle_{w_{h}(\Theta)} = \int^{2}_{0}\rho(\Theta) \varphi_{h,s}(\Theta)w_{h}(\Theta)d\Theta.
\end{align}
In equation (\ref{5}), $ \langle.,.\rangle $ denotes the inner product.\\
The truncated Lucas wavelet series can be defined as
\begin{align}\label{6}
\rho(\Theta)= \sum_{h=0}^{2^{k}-1}\sum_{s=0}^{S-1}\mathcal{E}_{h,s}\varphi_{h,s}(\Theta)= E^{T}\Psi(\Theta).
\end{align}

\noindent Here, $ \mathcal{E} $ and $ \Psi(\Theta) $ are $ 2^{k}S\times 1 $ matrices.
\begin{thm} \label{theorem1}
If $ \Psi(\Theta) $ is Lucas wavelets vector then,
\begin{align}\label{62b}
\dfrac{d\Psi(\Theta)}{d\Theta} = D\Psi(\Theta),
\end{align}
be the the $ I^{st} $ derivative of the vector $ \Psi(\Theta)$, where, $ D $ is $ 2^{k}S \times 2^{k}S$ 
 matrix which is defined by 
\begin{equation}
   D=
  \begin{pmatrix}
     F & \bf0 \\
  \bf0 & F \\
  \end{pmatrix}
  \label{62a}
\end{equation}
Also, $ F $ is a square matrix of order $ S $ whose $ (\Upsilon, \lambda)^{th} $ element is defined by 
\begin{align}\label{u1}
F_{\Upsilon, \lambda}=
\begin{cases}
2^{k+1}\iota (\Upsilon-1)\sqrt{\frac{\alpha_{\Upsilon-1}}{\alpha_{\lambda-1}}} ,\quad \Upsilon=2,\ldots,S,\,\,\, \lambda=1,\ldots,\Upsilon-1 \,\,\,{\rm and} \,\,\,(\Upsilon+\lambda)odd, \\  
  0,\hspace{2.92cm} o.w. 
  \end{cases}
  \end{align}
\end{thm}
  
\begin{proof}
The  proof of this theorem can be obtain by using previous definitions. 
\end{proof}

\begin{cor}
The differential operational matrix for $h^{th}$-order can be attained by utilizing (\ref{62b}) as
\begin{align}\label{c1}
\dfrac{d^{(h)}\Psi(\Theta)}{d\Theta^{(h)}} = D^{(h)}\Psi(\Theta),
\end{align}
where $ D^{(h)} $ is $ h^{th} $-order of the matrix $ D $.
\end{cor}

Further, we will introduce product operational matrix, which is utilized for solving the differential equations 
\begin{align}\label{m6}
\Psi(\Theta)\Psi^{T}(\Theta)E\cong \tilde{E}\Psi(\Theta).
\end{align}
We consider the second order differential equation for the foundation of shifted Lucas wavelet technique as\\
\begin{align}\label{s1}
\rho^{''}(\Theta)= G(\Theta,\rho(\Theta),\rho(\alpha \Theta)), \hspace{0.3cm} \Theta\in (0,l),
\end{align}
with supplementary conditions
\begin{align}\label{s2}
\rho(\Theta)\vert_{\Theta=0}=A_{1},\hspace{0.1cm}
\rho^{'}(\Theta)\vert_{\Theta=0}=A_{2},
\end{align}
\begin{align}\label{s3}
\rho(\Theta)\vert_{\Theta=0}=B_{1},\hspace{0.1cm}
\rho^{'}(\Theta)\vert_{\Theta=l}=B_{2}.
\end{align}
In order to apply the shifted orthogonal Lucas wavelet technique, firstly we have to approximate the unknown function $ \rho(\Theta) $ and $ \rho^{'}(\Theta) $, $ \rho^{''}(\Theta) $ in equation (\ref{s1}), (\ref{s2}) and (\ref{s3} and construct a residual function as: 
\begin{align}\label{e6}
\Re(\Theta)=  E^{T}D^{2}\Psi(\Theta)-G\left(\Theta,E^{T}\Psi(\Theta),E^{T}\Psi(\alpha\Theta) \right) .
\end{align} 
Considering the Tau method, equation (\ref{e6}) generates the
 $( 2^{k}S-2) $ nonlinear equations with the unknown coefficients $\mathcal{E} _{n,m} $.
The supplementary conditions are approximated as: 
\begin{align}\label{e8}
\rho(\Theta)\vert_{\Theta=0}=E^{T}\Psi(0)= A_{1}, \hspace{0.1cm} \rho^{'}(\Theta)\vert_{\Theta=0} = E^{T}D\Psi(0)= A_{2},
\end{align}
\begin{align}\label{e81}
\rho(\Theta)\vert_{\Theta=0}=E^{T}\Psi(0)= B_{1}, \hspace{0.1cm} \rho^{'}(\Theta)\vert_{\Theta=l} = E^{T}D\Psi(0)= B_{2},
\end{align}
Thus, approximated initial and boundary conditions and $( 2^{k}S-2) $ nonlinear equation generates set of $ 2^{k}S $ nonlinear equations. Further, these equations can be solved to determine the unknown component of vector $ E $. Hence $ \rho(\Theta) $ can be calculated.
\section{Convergence and error analysis}
\begin{thm}\label{th1}
The Lucas wavelet series solution defined in equation (\ref{4}) converges towards $\rho(\Theta)$ by using Lucas wavelet.
\end{thm} 
\begin{proof}
proof is omitted.
\end{proof}
\begin{thm}\label{th2}
Suppose that $ \rho(\Theta)\in L_{w_{h}}^{2}[0,2] $, as well as bounded second order derivative, say $ \left |\rho''(\Theta)\right |\leqslant N $, then the function $\rho(\Theta)$ can be elaborated as an infinite sum of Lucas wavelets and the series is uniformly convergent to $\rho(\Theta)$. Additionally, the following inequalities are satisfied by the coefficients $\mathcal{E}_{h,s}$ defined in \eqref{5}:
\begin{align}
\left\vert \mathcal{E}_{h,s}\right\vert&\leq\frac{2N\sqrt{\pi}}{(h+1)^{5/2}(s^{2}-1)},\quad \forall s>1,\,\,\,h\geq 0,\\
{\rm and} \quad \left\vert \mathcal{E}_{h,1}\right\vert&\leq\frac{\sqrt{\pi}}{(h+1)^{3/2}}\max\limits_{0\leq \Theta\leq 2}\left\vert \rho'(\Theta)\right\vert,\quad {\rm for} \,\,\,s=1, \forall h\geq 0.
\end{align}
\end{thm}
\begin{proof}
Proof is omitted.
\end{proof}
\begin{lem}\label{l1}
Let $ \Theta\geqslant s$, then $ f(\Theta) $ will be a continuous, positive and decreasing function. If $ f(k) = \mathcal{E}_{k} $ along with $ \sum \mathcal{E}_{h} $, is convergent and $ \Re_{h} = \sum\limits_{k=h+1}^{\infty}\mathcal{E}_{k}, $ then $ \Re_{h}\leq \int\limits_{h}^{\infty}f(\Theta)d\Theta$ \cite{Stewart2012}.  
\end{lem}
\begin{thm}
If $ \rho\in L^{2}_{w_{h}} $ defined in (\ref{4})  satisfies theorem (\ref{th1}) and 
\begin{align*}
\rho_{k,S}(\Theta)= \sum_{h=0}^{2^{k}-1}\sum_{s=0}^{S-1}\mathcal{E}_{h,s}\varphi_{h,s}(\Theta)
\end{align*}
then 
\begin{align*}
\parallel \rho-\rho_{k,S}\parallel_{w_{h}}< \sqrt{N^{2}\pi\left(\dfrac{1}{2^{5(2^{k}-1)-2}5ln (2)} \right) \left(\dfrac{(S^{2}-2S)ln(S)-S^{2}ln(S-2)+(2ln(S-2)-2)S+2}{4S(S-2)} \right) }
\end{align*}
is the error estimation, for $ S>2. $
\end{thm}
\begin{proof}
Proof is omitted.
\end{proof}

\section{Illustrative examples}
In this section, we consider the following test problems to show the validity and accuracy of the proposed scheme. \\
\textbf{Problem 1.} Let us consider the nonlinear Lane-Emden problem \cite{Aminikhah2013}
\begin{align}\label{e161}
\rho^{''}(\Theta)+\dfrac{6}{\Theta}\rho^{'}(\Theta)+14\rho(\Theta)=-4\rho(\Theta)ln(\rho(\Theta)), \hspace{0.2cm} 0< \Theta \leq 1,
\end{align}
with
\begin{align}\label{e162}
\rho(\Theta)\vert_{\Theta=0}=1, \hspace{0.2cm}\rho^{'}(\Theta)\vert_{\Theta=0}=0.
\end{align}
This problem has exact solution $ \rho(\Theta) = e^{-\Theta^{2}}. $ First of all to solve equation (\ref{e161}), we are used same transformation as we have applied in above example i.e. $ \rho(\Theta) $ = $ e^{z(\Theta)} $.  Then transformed form of equations (\ref{e161}) and (\ref{e162}) are 
\begin{align}\label{e163}
z^{''}(\Theta)+(z^{'}(\Theta))^{2}+\dfrac{6}{\Theta}z^{'}(\Theta)+14 = -4 z(\Theta), \hspace{0.1cm} \Theta\geq 0
\end{align}
with
\begin{align*}
z(\Theta)\vert_{\Theta=0}=0, \hspace{0.2cm} z^{'}(\Theta)\vert_{\Theta=0}=0.
\end{align*}
Also,
\begin{align*}
z(\Theta)=-\Theta^{2}.
\end{align*}
The values of component vector $ E $ for equation (\ref{e163}) are formed by utilizing  the method described in section (\ref{g}). In this case, we take $ k=0 $ and $ S=3 $. After some manipulations we get
\begin{align*}
\rho(\Theta) =e^{z(\Theta)}= e^{5.55112\times10^{-17}-\Theta^{2}},
\end{align*}
which is the approximate result for considered problem.\\
\textbf{Problem 2.} Let us consider the linear Pantograph initial value problem \cite{Rahimkhani2017}
\begin{align}\label{e16}
\rho^{''}(\Theta)-\dfrac{3}{4}\rho(\Theta)-\rho\left( \dfrac{\Theta}{2}\right) +\Theta^{2}-2=0, \hspace{0.2cm} 0\leq \Theta\leq 1,
\end{align}
with
\begin{align}
\rho(\Theta)\vert_{\Theta=0}=\rho^{'}(\Theta)\vert_{\Theta=0}=0.
\end{align}
This problem has exact solution $ \rho(\Theta) =\Theta^{2}$ at $ \gamma=2 $ given in \cite{Rahimkhani2017}. In the proposed technique let $ k=0 $ and $ S=3 $. After performing some calculations we get 
\begin{align*}
\rho(\Theta)=&E^{T}\Psi(\Theta)=(-4.996003610813204\times 10^{-16} + 8.088790061448873\times 10^{-17}\iota)&\\&+(2.220446049250313\times 10^{-16} + 
     2.7193641226633042\times 10^{-17}\iota) \Theta &\\&
     + (0.9999999999999994 - 
    4.8544160607900805\times 10^{-16} \iota) \Theta^2.
\end{align*}
which is the approximate solution for considered problem. 
\section{Conclusion}
In the present article, we have developed a novel technique for numerical solutions of differential equations. The basic foundations for this technique are shifted Lucas polynomials and Tau method. The shifted Lucas polynomials are well behaved orthonormal basis functions. To show the validity of NSLW technique, we have considered two different kind of differential equations. We hope that, the initiative of introducing the shifted Lucas polynomials will fill the gap of literature, related to Lucas polynomials. The computations associated with all the considered examples are done by software Mathematica 9.

\end{document}